\documentclass[12pt]{article}
\usepackage[utf8]{inputenc}
\usepackage[T1]{fontenc}
\usepackage{graphicx}
\usepackage{grffile}
\usepackage{longtable}
\usepackage{wrapfig}
\usepackage{rotating}
\usepackage[normalem]{ulem}
\usepackage{amsmath}
\usepackage{textcomp}
\usepackage{amssymb}
\usepackage{capt-of}
\usepackage[colorlinks=true]{hyperref}
\usepackage{minted}
\usepackage[a4paper,margin=1.5cm,includefoot]{geometry}
\usepackage{amsmath,amssymb,amsthm}
\usepackage{graphicx}
\newtheorem{theorem}{Theorem}[section]
\newtheorem{proposition}[theorem]{Proposition}
\newtheorem{corollary}[theorem]{Corollary}
\newtheorem{lemma}[theorem]{Lemma}
\theoremstyle{definition}
\newtheorem{definition}[theorem]{Definition}
\theoremstyle{remark}

\newcommand{\Haus}{\mathcal H}
\newcommand{\dimH}{\dim_{\mathrm H}}
\newcommand{\dimP}{\dim_{\mathrm P}}
\newcommand{\dimB}{\dim_{\mathrm B}}
\newcommand{\udimB}{\overline{\dim}_{\mathrm B}}
\newcommand{\ldimB}{\underline{\dim}_{\mathrm B}}
\newcommand{\lowM}{\underline{\mathcal M}}
\newcommand{\uppM}{\overline{\mathcal M}}
\newcommand{\eps}{\epsilon}
\newcommand{\cov}{\mathcal N}
\newcommand{\cost}{\mathcal C}
\author{Siniša Miličić}
\date{}
\title{Scaling cores and dust-free layered sets}
\hypersetup{
 pdfauthor={Siniša Miličić},
 pdftitle={Scaling cores and dust-free layered sets},
 pdfkeywords={},
 pdfsubject={},
 pdfcreator={Emacs 30.2 (Org mode 9.7.30)}, 
 pdflang={English}}
\begin{document}

\maketitle
\begin{center}
\small Faculty of Informatics, Juraj Dobrila University of Pula\\
Ulica Alda Negrija 6, 52100 Pula, Croatia\\
\texttt{smilicic@unipu.hr}\\
ORCID: \texttt{0000-0002-5659-5598}
\end{center}

\begin{abstract}
We study Hausdorff–box dimension gaps in compact sets of finite positive Hausdorff measure. At each scale, the points of low local measure form a scaling core that carries all covering growth above the Hausdorff dimension. This gives gap criteria for both box and intermediate dimensions. We construct dust-free layered sets, in which every neighbourhood of every point has positive measure. Their dimensions are determined by the layer dimension and a product core formed from a horizontal limit set and a sequence of heights. Under a suitable approximation-rate condition, the volume remaining after subtraction of the product tube grows linearly with radius, with coefficient twice the total layer measure. Every compact null horizontal set admits such a construction. Cantor approximants, homogeneous cube constructions and Bedford–McMullen carpets illustrate the results.
\end{abstract}

\setcounter{tocdepth}{1}
\tableofcontents
\section{The model and the dimension-gap problem}
\label{sec:orged263e7}

Finite positive Hausdorff measure does not force equality of Hausdorff and box dimensions. Put \(V_\alpha=\{0\}\cup\{n^{-\alpha}:n\ge1\}\) for \(\alpha>0\). When \(0<\alpha<1\), the compact set
\[
([0,1]\times\{0\})\cup(V_\alpha\times V_\alpha)
\]
is countably \(1\)-rectifiable, has length \(1\), Hausdorff and packing dimensions \(1\), and box dimension \(2/(\alpha+1)>1\). Countable dust overtakes the unit interval in box counting but carries no length \cite{mendivil2011}.

We call a compact \(F\) with \(0<\Haus^d(F)<\infty\) \emph{dust-free} if
\[
F=\operatorname{supp}(\Haus^d\!\restriction F).
\]
Thus every ball centred on \(F\) has positive \(d\)-dimensional measure; no uniform lower density bound is imposed. The question is whether a dimension gap can persist under this support condition, and how to identify the part responsible for it.

\begin{figure}[tbp]
\centering
\includegraphics[width=0.97\linewidth]{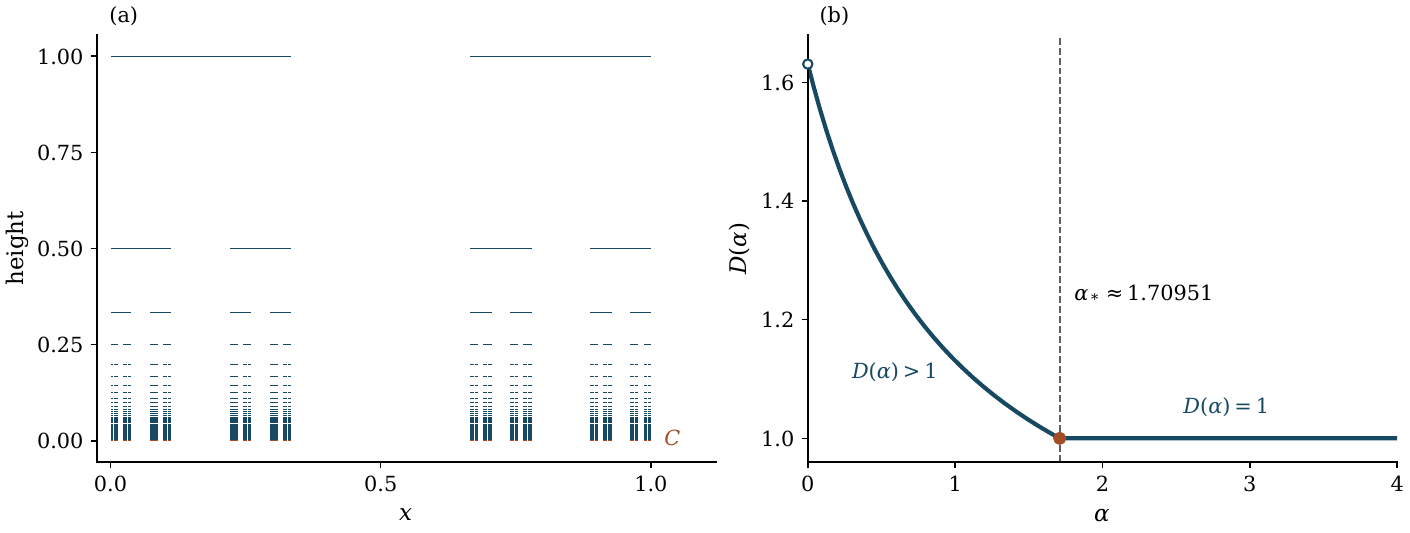}
\caption{Left: the closed grid squares of side \(3^{-6}\) meeting the full set \(F_1\), including every Cantor layer and the base (orange). Right: the box dimension of \(F_\alpha\), with transition at \(\alpha_*=(\log_3 2)/(1-\log_3 2)\).}
\label{fig:cantor}
\end{figure}

Our model places the middle-third Cantor approximants \(C_n\) at polynomial heights. We take \(C_0=[0,1]\), so \(|C_n|=(2/3)^n\), and start the layers at \(n=1\):
\begin{equation}\label{eq:cantor-model}
F_\alpha=\bigcup_{n\ge1}C_n\times\{n^{-\alpha}\}\ \cup\ (C\times\{0\}),
\qquad C=\bigcap_{n\ge1}C_n,\quad \alpha>0.
\end{equation}
\begin{equation}\label{eq:cantor-overview}
\Haus^1(F_\alpha)=2,\qquad
\dimH F_\alpha=\dimP F_\alpha=1,\qquad
\dimB F_\alpha=\max\left\{1,\log_3 2+\frac1{\alpha+1}\right\}.
\end{equation}
Although the base \(C\times\{0\}\) has zero length, every neighbourhood of a base point meets a positive-length piece of a sufficiently low layer. Hence \(F_\alpha\) is dust-free. Figure \ref{fig:cantor} shows the construction and its dimension transition; Corollary \ref{cor:power} proves \eqref{eq:cantor-overview}.

The common principle is that a remainder costing \(O(\eps^{-d})\) covering balls cannot affect dimensions above \(d\). For an arbitrary compact set of finite positive \(d\)-dimensional Hausdorff measure, Theorem \ref{thm:density-core} identifies the dimension gap through its scaling cores. For layered sets, Theorem \ref{thm:summable-approximation} gives a fixed product core; Lemma \ref{lem:layered-profiles} supplies covering comparisons uniform in the heights. The product-core dimensions determine the scaling-core exponents above \(d\). Theorem \ref{thm:additive-tube} recovers the layer measure from an additive tube formula with absolute error \(o(\eps)\), even when the product tube dominates the linear term. Corollary \ref{cor:null-core-realization} realizes every compact null horizontal set in this way.

The support condition alone does not control tube growth; Example 3 of Ambrosio, Colesanti and Villa \cite{ambrosio_colesanti_villa2008} also exhibits this issue. Tricot's nucleus–tail method \cite{tricot1995} is a foundation of the Žubrinić–Županović programme and subsequent developments on spirals, chirps and orbits \cite{zubrinic_zupanovic2005,elezovic_zupanovic_zubrinic2007,resman2013,rolin_vlah_zupanovic2021}. There, separated early pieces and a merged geometric nucleus may contribute at the same order. Our cores play an analogous role in locating the geometry responsible for tube growth: the scaling core selects points of low local measure, while the product core is the repeated horizontal limit set. The additive theorem separates the entire product tube from the contribution of the layers; the definitions differ from a nucleus formed by merging neighbouring pieces.

Generalized nests couple size and spacing \cite{milicic2018}, whereas here the horizontal approximations and vertical heights can be chosen independently. Related inhomogeneous attractors have dimension gaps with infinite layer measure \cite{fraser2012,fraser2016carpets}. We retain finite positive layer measure and obtain explicit formulas from sequence asymptotics and intermediate-dimension theory.
\section{Scaling cores and dimension gaps}
\label{sec:org2980450}

We first leave the layered model and consider arbitrary compact sets of finite positive \(d\)-dimensional Hausdorff measure. Points with sufficiently high density can be covered by \(O(\eps^{-d})\) balls. The remaining low-density points form the scaling core: they carry any covering growth above dimension \(d\).

The open tube set, or Minkowski sausage, of \(E\) is \((E)_\eps=\{x:\operatorname{dist}(x,E)<\eps\}\). Write \(B(x,r)\) for the open ball and \(|U|\) for ambient Lebesgue measure. Hausdorff measure is normalized to agree with volume on integer-dimensional planes. We use \(\dimH,\dimP,\ldimB,\udimB\) for Hausdorff, packing and lower/upper box dimensions \cite{falconer2014,mattila1995}.

For a bounded set \(E\), let \(\cov_\eps(E)\) be the least number of closed radius-\(\eps\) covering balls, with \(\cov_\eps(\varnothing)=0\). Standard packing and product arguments give
\begin{equation}\label{eq:covering-tools}
|(E)_\eps|\asymp\eps^m\cov_\eps(E)\quad(E\subset\mathbb R^m),
\qquad
\cov_\eps(X\times V)\asymp\cov_\eps(X)\cov_\eps(V).
\end{equation}
These comparisons hold at every positive scale with dimensional constants, uniformly in the bounded sets. Also \(\cov_{c\eps}(E)\asymp\cov_\eps(E)\) for each fixed \(c>0\). Tubes and covering numbers are exactly unchanged by closure: distance to a set equals distance to its closure, and every finite union of closed covering balls contains the closure of the set it covers.

Throughout, \(u_\eps\asymp v_\eps\) means \(c v_\eps\le u_\eps\le C v_\eps\) for some \(0<c\le C<\infty\) and all sufficiently small \(\eps>0\); \(u\sim v\) means that the ratio tends to one. Constants and thresholds may depend on fixed data but not on \(\eps\); any further uniformity is specified. Limits in other variables are stated explicitly. The same convention applies to \(O(\cdot)\).

For \(0<\eps<1\), \(0<\theta\le1\), and \(t\ge0\), define the admissible covering cost
\begin{equation}\label{eq:cover-cost}
\cost_{\eps,\theta}^t(E)
=\inf\left\{\sum_i(\operatorname{diam}U_i)^t:
E\subset\bigcup_iU_i,\quad
\eps\le\operatorname{diam}U_i\le\eps^\theta\right\}.
\end{equation}
The empty set has cost zero. For fixed \(0<\theta<1\), both endpoints tend to zero, but a range of diameters is allowed. At \(\theta=1\), every covering set has diameter exactly \(\eps\), and
\begin{equation}\label{eq:cost-tube-compatibility}
\cost_{\eps,1}^t(E)\asymp\eps^t\cov_\eps(E)
\asymp\eps^{t-m}|(E)_\eps|.
\end{equation}
The constants depend only on \(m\), uniformly in \(E,t\) and \(0<\eps<1\): a diameter-\(\eps\) set lies in a radius-\(\eps\) ball, and each such ball is covered by a dimensional constant number of diameter-\(\eps\) sets.

The lower \(\theta\)-intermediate dimension is the infimum of \(t\) for which \(\cost_{\eps,\theta}^t(E)\) has lower limit zero as \(\eps\downarrow0\); the upper version uses upper limit zero \cite{falconer2020}. Thus \(\theta=1\) gives the box dimensions. We write \(\dim_\theta\) or \(\dimB\) when the corresponding lower and upper values agree. We define \(\dim_0=\dimH\) separately; we do not substitute \(\theta=0\) into the shrinking-scale argument, nor assume continuity at zero.

\begin{definition}[Scaling core]\label{def:density-core}
Let \(A\subset\mathbb R^m\) be compact, \(0<d<m\), and \(0<M:=\Haus^d(A)<\infty\). Put \(\mu=\Haus^d\!\restriction A\). For fixed \(\kappa>0\), its scaling core at scale \(\eps\) is
\begin{equation}\label{eq:density-core}
K_\eps=K_\eps(A;\kappa)
=\{x\in A:\mu(B(x,\eps))\le\kappa\eps^d\}.
\end{equation}
For \(0<\theta\le1\), define its moving-family exponents by
\begin{align}
\underline D_{\theta,\kappa}(A)
&=\inf\{t\ge0:\liminf_{\eps\downarrow0}\cost_{\eps,\theta}^t(K_\eps)=0\},\label{eq:core-exponents}\\
\overline D_{\theta,\kappa}(A)
&=\inf\{t\ge0:\limsup_{\eps\downarrow0}\cost_{\eps,\theta}^t(K_\eps)=0\}.\nonumber
\end{align}
\end{definition}

The core is compact: Fatou's lemma makes \(x\mapsto\mu(B(x,\eps))\) lower semicontinuous. The core-defining radius is the smallest permitted covering diameter. These exponents describe a moving family, not a fixed limiting set. By \eqref{eq:cost-tube-compatibility}, at \(\theta=1\) they are the lower and upper limits of \(\log\max\{1,\cov_\eps(K_\eps)\}/\log(1/\eps)\).

The next lemma separates the geometric estimates from the dimension criterion. The complement has enough mass in every radius-\(\eps\) ball to limit the number of disjoint such balls.

\begin{lemma}[Covering the complement]\label{lem:scaling-complement}
Under Definition~\ref{def:density-core}, for every \(\eps>0\), the complement \(A\setminus K_\eps\) is coverable by at most \(C_m\kappa^{-1}M\eps^{-d}\) sets of diameter exactly \(\eps\), and
\begin{equation}\label{eq:density-core-tube}
0\le |(A)_\eps|-|(K_\eps)_\eps|
\le C_m\kappa^{-1}M\eps^{m-d}.
\end{equation}
For every \(0<\eps<1\), \(0<\theta\le1\), and \(t\ge0\),
\begin{equation}\label{eq:density-core-cost}
0\le\cost_{\eps,\theta}^t(A)-\cost_{\eps,\theta}^t(K_\eps)
\le C_m\kappa^{-1}M\eps^{t-d}.
\end{equation}
All constants depend only on the ambient dimension \(m\).
\end{lemma}
\begin{proof}
Let \(P_\eps\) be a maximal strictly \(2\eps\)-separated subset of \(A\setminus K_\eps\); boundedness makes it finite. Its open radius-\(\eps\) balls are disjoint. Since each centre lies outside the core,
\[
\#P_\eps\,\kappa\eps^d
\le\sum_{p\in P_\eps}\mu(B(p,\eps))
\le\mu(A)=M.
\]
Maximality gives a cover of the complement by the closed balls \(\overline B(p,2\eps)\). A grid of mesh \(\eps/\sqrt m\) covers each such ball with a dimensional constant number of cubes of diameter \(\eps\). This proves the covering claim.

Every point of \((A)_\eps\setminus(K_\eps)_\eps\) is within \(\eps\) of a point outside the core, hence within \(3\eps\) of some \(p\in P_\eps\). Thus the tube excess lies in \(\bigcup_{p\in P_\eps}B(p,3\eps)\), giving \eqref{eq:density-core-tube}. For \eqref{eq:density-core-cost}, append the diameter-\(\eps\) cover of the complement to an admissible cover of \(K_\eps\). The additional \(t\)-cost is at most \(C_m\kappa^{-1}M\eps^{t-d}\); take the infimum. Both lower bounds follow from \(K_\eps\subset A\).
\end{proof}

\begin{theorem}[Scaling-core criterion]\label{thm:density-core}
Under Definition~\ref{def:density-core}, fix \(0<\theta\le1\). The lower or upper \(\theta\)-intermediate dimension of \(A\) exceeds \(\dimH A=d\) if and only if the corresponding scaling-core exponent exceeds \(d\). More precisely,
\begin{align}
\underline{\dim}_\theta A&=\max\{d,\underline D_{\theta,\kappa}(A)\},\label{eq:core-lower-dimension}\\
\overline{\dim}_\theta A&=\max\{d,\overline D_{\theta,\kappa}(A)\}.
\label{eq:core-upper-dimension}
\end{align}
The covering numbers satisfy
\begin{equation}\label{eq:density-core-count}
\cov_\eps(A)\asymp\eps^{-d}+\cov_\eps(K_\eps)
\qquad(\eps\downarrow0),
\end{equation}
with constants and threshold depending on \(A,d,m,\kappa\).
\end{theorem}
\begin{proof}
Lemma~\ref{lem:scaling-complement} gives
\(\cov_\eps(A)\le\cov_\eps(K_\eps)+C_m\kappa^{-1}M\eps^{-d}\).
For the reverse bound, \(K_\eps\subset A\), while any radius-\(\eps\) ball cover of \(A\) is a Hausdorff cover with diameters at most \(2\eps\). The corresponding Hausdorff contents tend to \(M>0\); hence \(\cov_\eps(A)\eps^d\ge c_{A,d}>0\) for sufficiently small \(\eps\). Combining these two lower bounds proves \eqref{eq:density-core-count}.

For every \(t>d\), \eqref{eq:density-core-cost} bounds the difference of the two nonnegative covering costs by a quantity tending to zero. Their lower limits therefore vanish together, and their upper limits vanish together. The vanishing-cost thresholds coincide above \(d\). Since both intermediate dimensions of \(A\) are at least \(\dimH A=d\), the two maximum formulas follow.
\end{proof}

Dust-freeness is not needed here. For a dust-free set, every centred ball has positive mass, but that mass need not have a uniform lower bound relative to \(\eps^d\). Every fixed \(\kappa>0\) detects the same dimension gap. At \(t=d\), the cost error is only \(O(1)\), so the argument does not identify critical contents.

If \(|(A)_\eps|/\eps^{m-d}\to\infty\), the tube estimate gives \(|(K_\eps)_\eps|/|(A)_\eps|\to1\). A lower-box gap guarantees this throughout the limit; an upper-box gap guarantees it along suitable sequences.
\section{Summable layered sets}
\label{sec:orgb027547}

We now approximate a compact horizontal limit set by sets placed at decreasing heights. The integer \(d\ge1\) is the horizontal dimension and will be the Hausdorff dimension of the layered set; the approximated set may have smaller dimension. The geometry supplies a fixed product core. Its controlled remainder is verified directly, since the scaling-core criterion concerns a different, measure-defined family.

\begin{definition}[Summable approximation scheme]\label{def:summable-approximation}
Let \(X\subset\mathbb R^d\) be nonempty, compact and Lebesgue-null. A summable approximation scheme consists of compact sets \(X_n\) and positive \(\delta_n\downarrow0\) such that
\begin{align}
X&\subset X_n\subset\overline{(X)_{\delta_n}},\label{eq:approximation-sandwich}\\
S_X:=\sum_{n\ge1}|(X)_{2\delta_n}|&<\infty,
\qquad S:=\sum_{n\ge1}|X_n|>0.
\label{eq:neighborhood-summability}
\end{align}
Choose heights \(h_1>h_2>\cdots>0\) tending to zero. Define the height set, layered set and fixed product core by
\begin{align}
V&=\{0,h_1,h_2,\ldots\},\nonumber\\
F&=\bigcup_{n\ge1}X_n\times\{h_n\}\ \cup\ (X\times\{0\}),\label{eq:layered-set}\\
G&=X\times V
=\overline{\bigcup_{n\ge1}X\times\{h_n\}}.\label{eq:product-core}
\end{align}
We call \(F\) the associated summable layered set.
\end{definition}

Each \(X_n\) contains \(X\) and stays within distance \(\delta_n\) of it. The sum \(S_X\) requires the volumes of these enclosing tubes to be summable; the necessary rate of \(\delta_n\) depends on the tube decay of \(X\). It implies \(S\le S_X<\infty\), while \(S>0\) requires only one positive-volume approximant. The product core \(G\) is a spine consisting of a copy of \(X\) at every height in \(V\), including zero.

The following lemma isolates the transfer argument for a fixed core and a remainder of controlled covering size. It allows the retained part to lie near the core rather than inside it.

\begin{lemma}[Removal of a controlled remainder]\label{lem:cover-transfer}
Let \(0<d\le m\), let \(E,H\subset\mathbb R^m\) be bounded, with \(0<\Haus^d(E)<\infty\) and \(H\subset\overline E\). Suppose there exist \(C>0\) and, for every sufficiently small \(\eta>0\), a set \(R_\eta\) such that
\[
E\subset R_\eta\cup(H)_{C\eta},\qquad
\cov_\eta(R_\eta)\le C\eta^{-d}.
\]
Then \(\cov_\eps(E)\asymp\eps^{-d}+\cov_\eps(H)\) as \(\eps\downarrow0\), and, for \(0<\theta\le1\),
\[
\underline{\dim}_\theta E=\max\{d,\underline{\dim}_\theta H\},\qquad
\overline{\dim}_\theta E=\max\{d,\overline{\dim}_\theta H\}.
\]
\end{lemma}
\begin{proof}
The remainder bound and fixed-radius covering comparison give the covering upper bound. Positivity of \(\Haus^d(E)\) gives \(\cov_\eps(E)\ge c_E\eps^{-d}\) at small scales; \(H\subset\overline E\) gives the other lower bound.

For intermediate dimensions, enlarge each set in an admissible diameter-\(r_i\) cover of \(H\) by \(C\eta\), where \(r_i\in[\eta,\eta^\theta]\). Euclidean doubling covers this enlargement by a bounded number of sets of diameter \(r_i\). Cover \(R_\eta\) by \(O(\eta^{-d})\) sets of diameter \(\eta\). The resulting \(t\)-cost is at most
\(C'(\sum_i r_i^t+\eta^{t-d})\).
For \(t>d\), the error tends to zero, proving both upper bounds with their respective scale quantifiers. The lower bounds follow from \(\dimH E=d\), monotonicity and closure invariance for \(\theta>0\).
\end{proof}

\begin{lemma}[Layered covering estimates]\label{lem:layered-profiles}
For the scheme of Definition~\ref{def:summable-approximation}, as \(\eps\downarrow0\),
\begin{align}
\cov_\eps(F)&\asymp\eps^{-d}+\cov_\eps(X)\cov_\eps(V),\label{eq:approximation-count}\\
|(F)_\eps|&\asymp\eps+\eps^{d+1}\cov_\eps(X)\cov_\eps(V).\label{eq:approximation-tube}
\end{align}
The comparison constants and threshold depend only on the approximation scheme, uniformly in the heights.
\end{lemma}
\begin{proof}
For \(\delta_n\ge\eps\), take a maximal strictly \(\eps\)-separated subset of \(X_n\). Its disjoint radius-\(\eps/2\) balls lie in \((X)_{2\delta_n}\), so
\[
\cov_\eps(X_n)\le C_d\eps^{-d}|(X)_{2\delta_n}|.
\]
Thus the union \(R_\eps\) of these layers has a cover by at most \(C_dS_X\eps^{-d}\) balls, independently of their heights. Every remaining layer, and the base, lies in \((G)_\eps\). Covering this tube by balls of comparable radius and using \eqref{eq:covering-tools} proves the upper bound.

For the lower bound, \(G\subset F\). Also choose a fixed \(n_*\) with \(|X_{n_*}|>0\). Projecting any radius-\(\eps\) cover of that slice onto its horizontal plane gives
\[
\cov_\eps(F)\ge\cov_\eps(X_{n_*})
\ge\frac{|X_{n_*}|}{\omega_d\eps^d},\qquad \eps>0,
\]
where \(\omega_d\) is the volume of the unit ball in \(\mathbb R^d\). This bound is independent of all heights and holds at every scale. The tube comparison now follows from \eqref{eq:covering-tools} in ambient dimension \(d+1\).
\end{proof}

\begin{theorem}[Layered realization]\label{thm:summable-approximation}
The set \(F\) in Definition~\ref{def:summable-approximation} is compact, countably \(d\)-rectifiable, and
\begin{equation}\label{eq:approximation-measure}
\Haus^d(F)=S\in(0,\infty),\qquad \dimH F=\dimP F=d.
\end{equation}
For every \(0\le\theta\le1\),
\begin{align}
\underline{\dim}_\theta F&=\max\{d,\underline{\dim}_\theta G\},\label{eq:approximation-lower-intermediate}\\
\overline{\dim}_\theta F&=\max\{d,\overline{\dim}_\theta G\}.\label{eq:approximation-upper-intermediate}
\end{align}
\end{theorem}
\begin{proof}
The set is bounded. For a convergent sequence of its points, either a subsequence lies in one compact slice, or its layer indices tend to infinity. In the latter case the heights tend to zero and the distances of the horizontal coordinates from \(X\) tend to zero. The limit therefore lies in \(X\times\{0\}\), proving compactness.

The slices are disjoint, Hausdorff measure agrees with horizontal volume, and \(|X|=0\). Countable additivity gives \(\Haus^d(F)=\sum_n|X_n|=S\). A positive-volume slice gives dimension at least \(d\), while the countable cover by horizontal planes gives Hausdorff and packing dimensions at most \(d\), as well as countable \(d\)-rectifiability.

The remainder \(R_\eps\) constructed in Lemma~\ref{lem:layered-profiles} satisfies the hypotheses of Lemma~\ref{lem:cover-transfer} with \(E=F\) and \(H=G\). This proves the intermediate formulas for \(\theta>0\). At \(\theta=0\), countability of \(V\) gives \(\dimH G=\dimH X\le d\).
\end{proof}

\begin{corollary}[Dust-free realization]\label{cor:dust-free}
If every \(X_n\) is the support of its restricted Lebesgue measure, then \(F\) is dust-free.
\end{corollary}
\begin{proof}
Every ball centred on a positive-height slice meets a positive-volume part of that slice. For \((x,0)\in X\times\{0\}\) and \(r>0\), choose \(n\) with \(h_n<r/2\). Since \(x\in X_n\), full support gives \(|X_n\cap B(x,r/2)|>0\); its lift to height \(h_n\) lies in the radius-\(r\) ball about \((x,0)\).
\end{proof}

Combining Theorems \ref{thm:density-core} and \ref{thm:summable-approximation} relates the moving scaling cores to the fixed product core: for every fixed \(\kappa>0\) and \(0<\theta\le1\),
\begin{align}
\max\{d,\underline D_{\theta,\kappa}(F)\}
&=\max\{d,\underline{\dim}_\theta G\},\label{eq:density-product-bridge}\\
\max\{d,\overline D_{\theta,\kappa}(F)\}
&=\max\{d,\overline{\dim}_\theta G\}.\nonumber
\end{align}
Whenever the corresponding dimension of \(G\) exceeds \(d\), it equals the scaling-core exponent. In particular,
\begin{equation}\label{eq:product-core-gap}
\ldimB F>d\iff\ldimB G>d,
\qquad
\udimB F>d\iff\udimB G>d.
\end{equation}
The threshold is \(d=\dimH F\), not \(\dimH G\). For irregular factors the lower box dimension of \(G\) need not be the sum of the two lower box dimensions; \eqref{eq:approximation-count} keeps their scale dependence.
\section{The additive tube formula}
\label{sec:org3ebb7f9}

Fix the approximation scheme and heights of Definition \ref{def:summable-approximation}. We retain the layered set \(F\subset\mathbb R^{d+1}\), its product core \(G=X\times V\), and its total layer measure \(S=\sum_n|X_n|\). Thus \(G\) includes a copy of \(X\) at every height, not only at the base. The factor \(2\eps\) below comes from the one direction transverse to each layer.

\begin{lemma}[A compact flat layer]\label{lem:flat-slice}
Let \(Y\subset\mathbb R^d\) be compact and \(H=Y\times\{h\}\). Then
\begin{equation}\label{eq:flat-slice-tube}
2\eps|Y|\le|(H)_\eps|\le2\eps|(Y)_\eps|.
\end{equation}
For every continuous \(\varphi\) on a compact neighbourhood of \(H\),
\[
\frac1{2\eps}\int_{(H)_\eps}\varphi(z)\,dz
\longrightarrow\int_Y\varphi(x,h)\,dx.
\]
\end{lemma}
\begin{proof}
The cylinder \(Q_\eps=Y\times(h-\eps,h+\eps)\) satisfies
\[
Q_\eps\subset(H)_\eps\subset(Y)_\eps\times(h-\eps,h+\eps).
\]
This proves \eqref{eq:flat-slice-tube}. Compactness gives \(|(Y)_\eps\setminus Y|\to0\), hence
\[
|(H)_\eps\setminus Q_\eps|
\le2\eps\bigl(|(Y)_\eps|-|Y|\bigr)=o(\eps).
\]
The normalized integrals over the tube and cylinder therefore differ by \(o(1)\). Fubini's theorem and uniform continuity of \(\varphi\) give
\[
\frac1{2\eps}\int_{Q_\eps}\varphi(z)\,dz
=\int_Y\frac1{2\eps}\int_{h-\eps}^{h+\eps}\varphi(x,t)\,dt\,dx
\longrightarrow\int_Y\varphi(x,h)\,dx.
\]
No regularity of the boundary of \(Y\) is required.
\end{proof}

Infinitely many layer tubes may overlap, so the single-layer asymptotic cannot simply be summed. Fix \(L>2\). Summability controls the layers with \(\delta_n\ge\eps\), while those with \(\delta_n\le\eps^L\) fit into a negligible shell around the product tube. The next condition controls the intervening scales.

\begin{theorem}[Additive tube formula]\label{thm:additive-tube}
For a summable layered set, suppose that for some fixed \(L>2\),
\begin{equation}\label{eq:scale-window}
w_\eps|(X)_{2\eps}|\longrightarrow0,
\qquad w_\eps=\#\{n:\eps^L<\delta_n<\eps\},
\qquad\eps\downarrow0.
\end{equation}
Then
\begin{equation}\label{eq:additive-tube}
|(F)_\eps|=|(G)_\eps|+2\eps\Haus^d(F)+o(\eps).
\end{equation}
The limit is for fixed approximation scheme and heights; convergence is not asserted uniformly over height sequences.
\end{theorem}

Because \(G\subset F\), this recovers \(S=\Haus^d(F)\) from the excess over the product tube. For the Cantor model, the correction is \(4\eps+o(\eps)\), even when the product tube dominates it.

We use the shell estimate
\begin{equation}\label{eq:tube-shell}
|(K)_{r+u}|-|(K)_r|\le C_m(u/r)|(K)_r|,
\qquad0<u\le r,\quad K\subset\mathbb R^m\text{ compact}.
\end{equation}
For finite \(K\), contraction by \(r/t\) about each centre maps its radius-\(t\) tube inside its Voronoi cell into the radius-\(r\) tube in that cell. Summing gives \(|(K)_t|\le(t/r)^m|(K)_r|\), \(t\ge r\). Finite nets extend this to compact \(K\); subtraction gives \eqref{eq:tube-shell}. This is also a standard tube-volume estimate \cite{rataj_winter2010}.

\begin{proof}[Proof of Theorem~\ref{thm:additive-tube}]
Put \(E_\eps=(F)_\eps\setminus(G)_\eps\). Fix \(N\), set
\(J_N=\bigcup_{n=1}^N X_n\times\{h_n\}\), and put
\(E_{\eps,N}=(J_N)_\eps\setminus(G)_\eps\).
For small \(\eps\), the first \(N\) slice tubes are mutually disjoint. Their total overlap with \((G)_\eps\) is at most \(2N\eps|(X)_\eps|=o(\eps)\). Lemma~\ref{lem:flat-slice} therefore gives
\[
\frac{|E_{\eps,N}|}{2\eps}\longrightarrow
S_N:=\sum_{n=1}^N|X_n|.
\]
The remaining excess comes from layers with \(n>N\). Those with \(\delta_n\ge\eps\) contribute at most
\(2\eps\sum_{n>N}|(X)_{2\delta_n}|\).
Those with \(\eps^L<\delta_n<\eps\) contribute at most
\(2\eps w_\eps|(X)_{2\eps}|=o(\eps)\).
Layers with \(\delta_n\le\eps^L\) lie within distance \(\eps^L\) of \(G\); their excess is at most
\[
|(G)_{\eps+\eps^L}|-|(G)_\eps|
\le C\eps^{L-1}|(G)_\eps|=o(\eps),
\]
since \(L>2\) and the tube volumes are bounded. Consequently,
\begin{equation}\label{eq:tail-excess}
\limsup_{\eps\downarrow0}
\frac{|E_\eps\setminus E_{\eps,N}|}{2\eps}
\le T_N:=\sum_{n>N}|(X)_{2\delta_n}|.
\end{equation}
Thus \(S_N\le\liminf |E_\eps|/(2\eps)
\le\limsup |E_\eps|/(2\eps)\le S_N+T_N\).
Let \(N\to\infty\): then \(S_N\to S\) and \(T_N\to0\), proving \eqref{eq:additive-tube}.
\end{proof}

\begin{corollary}[Local tube excess]\label{cor:local-excess}
Under the hypotheses of Theorem~\ref{thm:additive-tube}, for every continuous \(\varphi\) on a fixed compact neighbourhood of \(F\),
\begin{equation}\label{eq:local-excess}
\frac1{2\eps}\int_{(F)_\eps\setminus(G)_\eps}\varphi(z)\,dz
\longrightarrow\int_F\varphi\,d\Haus^d.
\end{equation}
\end{corollary}
\begin{proof}
For fixed \(N\), Lemma~\ref{lem:flat-slice} and the \(o(\eps)\) overlap estimate in the preceding proof give
\[
\frac1{2\eps}\int_{E_{\eps,N}}\varphi(z)\,dz
\longrightarrow\int_{J_N}\varphi\,d\Haus^d.
\]
By \eqref{eq:tail-excess}, the omitted normalized integral has absolute upper limit at most \(\|\varphi\|_\infty T_N\). The omitted layer integral is at most \(\|\varphi\|_\infty(S-S_N)\) in absolute value. Both bounds tend to zero as \(N\to\infty\).
\end{proof}

Geometric approximation scales \(\delta_n=c\rho^n\), \(0<\rho<1\), satisfy \(w_\eps=O(\log(1/\eps))\); thus \eqref{eq:scale-window} holds whenever \(|(X)_r|=O(r^a)\) for some \(a>0\). Every compact null set also admits a sufficiently sparse scheme:

\begin{corollary}[Realization of an arbitrary null set]\label{cor:null-core-realization}
Let \(X\subset\mathbb R^d\) be nonempty, compact and null. For any fixed \(L>2\), choose positive scales satisfying
\[
\delta_1<1,\qquad \delta_{n+1}\le\delta_n^L,
\qquad |(X)_{2\delta_n}|\le2^{-n}.
\]
Then \(X_n=\overline{(X)_{\delta_n}}\) gives a dust-free summable layered set satisfying the layered realization theorem and the additive tube formula, for every strictly decreasing height sequence tending to zero.
\end{corollary}
\begin{proof}
Nullity permits the recursive choice. Each \(X_n\) is the closure of its interior and has positive volume. The sandwich and summability conditions follow immediately. If \(\delta_n<\eps\), then \(\delta_{n+1}<\eps^L\), so \(w_\eps\le1\) and \eqref{eq:scale-window} holds.
Apply Corollary~\ref{cor:dust-free} for full support.
\end{proof}

\begin{proposition}[Changing the bottom]\label{prop:bottom}
For bounded \(B\subset\mathbb R^d\), put
\(F_B=\bigcup_{n\ge1}X_n\times\{h_n\}\cup(B\times\{0\})\).
Then \(\Haus^d(F_B)=S+\Haus^d(B)\), using outer measure if necessary, and both Hausdorff and packing dimensions are \(d\). Its covering numbers are comparable to those of \(F\), and its lower and upper intermediate dimensions equal those of \(F\) for \(0\le\theta\le1\). Under \eqref{eq:scale-window},
\[
|(F_B)_\eps|=|(G)_\eps|+2\eps(S+|\overline B|)+o(\eps).
\]
\end{proposition}
\begin{proof}
Disjoint slices give the measure identity and the dimension statements at \(\theta=0\). A bounded bottom costs \(O(\eps^{-d})\) balls, while a positive-volume layer already requires that many. Lemma~\ref{lem:cover-transfer}, with core \(F\subset\overline{F_B}\), gives the intermediate identities. For the tube formula, use closure invariance and include \(\overline B\times\{0\}\) in the fixed finite family in the proof of Theorem~\ref{thm:additive-tube}. Its overlap with the product tube is at most \(2\eps|(X)_\eps|=o(\eps)\), and the same tail estimate applies.
\end{proof}

The distinction between \(B\) and \(\overline B\) is essential: a countable dense bottom has zero Hausdorff measure but the same tubes as its closure. This modification need not preserve dust-freeness.
\section{Prefractal applications}
\label{sec:orgf4db4d3}

We return to the canonical layered set \(F\) with base \(X\times\{0\}\). We first choose horizontal approximants \(X_n\) satisfying both layered theorems, then vary the height set \(V\). This gives dimension transitions and, when the product tube is negligible, exact Minkowski content.

For a bounded \(E\subset\mathbb R^m\), write
\[
\lowM^D(E)=\liminf_{\eps\downarrow0}\frac{|(E)_\eps|}{\eps^{m-D}},
\qquad
\uppM^D(E)=\limsup_{\eps\downarrow0}\frac{|(E)_\eps|}{\eps^{m-D}}.
\]
Positive finite lower and upper contents mean Minkowski nondegeneracy. If they coincide, \(E\) is Minkowski measurable, with unnormalized content \(\mathcal M^D(E)\).
\subsection{Cubes and Bedford–McMullen carpets}
\label{sec:org39138ab}

For cubes, choose integers \(d\ge1\), \(b\ge2\), retain \(1\le M<b^d\) closed \(b\)-adic subcubes of \([0,1]^d\), and iterate their homotheties \(S_i\). For a word \(w=(i_1,\ldots,i_n)\in\{1,\ldots,M\}^n\), write \(|w|=n\) and \(S_w=S_{i_1}\circ\cdots\circ S_{i_n}\). Put
\begin{equation}\label{eq:cubes}
X_n=\bigcup_{|w|=n}S_w([0,1]^d),\qquad
X=\bigcap_nX_n,\qquad s=\log_bM,\qquad q=M/b^d.
\end{equation}
The open set condition gives \(\dimH X=\dimP X=\dimB X=s\) and \(\cov_r(X)\asymp r^{-s}\) as \(r\downarrow0\) \cite{hutchinson1981}.

For Bedford–McMullen carpets, choose integers \(2\le b<c\), retain \(1\le M<bc\) rectangles in the \(b\)-column, \(c\)-row subdivision, and iterate their affine maps. Let \(X_n\) be the level-\(n\) union, \(X=\bigcap_nX_n\), \(d=2\), and \(M_1,\ldots,M_R\) the positive column counts. Set
\begin{equation}\label{eq:bm-exponents}
q=\frac M{bc},\qquad
s=s_B=\frac{\log R}{\log b}+\frac{\log(M/R)}{\log c}<2.
\end{equation}
Approximate squares give \(\dimP X=\dimB X=s_B\) and \(\cov_r(X)\asymp r^{-s_B}\), while
\[
\dimH X=\frac{\log\left(\sum_{i=1}^RM_i^{\log b/\log c}\right)}{\log b}.
\]
This equals \(s_B\) precisely when the occupied-column counts are equal \cite{mcmullen1984,fraser2021carpets}. All horizontal comparison constants depend only on the chosen construction.

In both schemes, each construction cell meets \(X\), has diameter at most \(c_0b^{-n}\), and has interior disjoint from the other level-\(n\) cells. Thus \(|X_n|=q^n\), \(X_n=\overline{\operatorname{int}X_n}\), and
\begin{equation}\label{eq:geometric-approximation}
X\subset X_n\subset\overline{(X)_{\delta_n}},
\qquad \delta_n=c_0b^{-n},
\end{equation}
where \(c_0=\sqrt d\) for cubes and \(c_0=2\) for carpets. Since \(|(X)_r|=O(r^{d-s})\) and \(s<d\), these neighborhoods are summable. There are only \(O(\log(1/\eps))\) approximation scales between \(\eps^L\) and \(\eps\), so the additional condition \eqref{eq:scale-window} also holds.

\begin{corollary}[Prefractal profiles]\label{cor:prefractal-profile}
For either scheme and heights \(1\ge h_1>h_2>\cdots\to0\), the resulting \(F\) is dust-free, with
\begin{align}
\Haus^d(F)&=\frac q{1-q},\qquad \dimH F=\dimP F=d,\label{eq:full-measure}\\
\cov_\eps(F)&\asymp\eps^{-d}+\eps^{-s}\cov_\eps(V),\label{eq:full-cover}\\
|(F)_\eps|&\asymp\eps+\eps^{d+1-s}\cov_\eps(V).\label{eq:full-tube}
\end{align}
These comparisons are uniform in the heights. In particular,
\[
\ldimB F=\max\{d,s+\ldimB V\},\qquad
\udimB F=\max\{d,s+\udimB V\}.
\]
For each fixed height sequence,
\begin{equation}\label{eq:prefractal-additive}
|(F)_\eps|=|(X\times V)_\eps|+\frac{2q}{1-q}\eps+o(\eps).
\end{equation}
Consequently,
\begin{equation}\label{eq:exact-content}
\eps^{d-s}\cov_\eps(V)\longrightarrow0
\quad\Longrightarrow\quad
\mathcal M^d(F)=\frac{2q}{1-q}.
\end{equation}
\end{corollary}
\begin{proof}
Apply Lemma~\ref{lem:layered-profiles}, Theorems~\ref{thm:summable-approximation} and~\ref{thm:additive-tube}, and Corollary~\ref{cor:dust-free}, using \(\sum_{n\ge1}q^n=q/(1-q)\) and the horizontal covering estimates. The last implication follows from \(|(X\times V)_\eps|/\eps\asymp\eps^{d-s}\cov_\eps(V)\).
\end{proof}
\subsection{Choosing the heights}
\label{sec:org802c185}

Power and logarithmic height laws are governed by the sequence and gap methods of \cite{elezovic_zupanovic_zubrinic2007,mardesic_resman_zupanovic2012,eichinger_winter2019}. For intermediate dimensions we use the decreasing-gap formula, Theorem 4.3 of Angelini and Molter \cite{angelini_molter2026}.

\begin{lemma}[Power-logarithmic heights]\label{lem:heights}
Fix \(\alpha>0\), \(\lambda\in\mathbb R\), and choose an integer \(n_0\ge2\) sufficiently large that
\(h_n=(n+n_0)^{-\alpha}(\log(n+n_0))^\lambda\) decreases strictly, with \(h_1\le1\). For \(\beta=(\alpha+1)^{-1}\),
\begin{align}
\cov_\eps(V)&\asymp\eps^{-\beta}(\log(1/\eps))^{\lambda\beta},
\qquad\eps\downarrow0,\label{eq:height-log}\\
\dim_\theta V&=\frac\theta{\alpha+\theta},\qquad0\le\theta\le1.
\label{eq:height-intermediate}
\end{align}
Constants depend on \(\alpha,\lambda,n_0\). For \(h_n=n^{-\alpha}\), the conclusions hold with \(\lambda=0\), without a shift.
\end{lemma}
\begin{proof}
The gaps satisfy \(h_n-h_{n+1}\sim\alpha n^{-\alpha-1}(\log n)^\lambda\) as \(n\to\infty\). Put \(N_\eps=\lceil\eps^{-\beta}(\log(1/\eps))^{\lambda\beta}\rceil\). For a sufficiently small fixed \(c>0\), this asymptotic bounds the gaps by more than \(2\eps\) for all sufficiently large indices up to \(\lfloor cN_\eps\rfloor\). Taking \(\eps\) smaller than half the minimum of the finitely many initial gaps gives separation throughout this initial segment. Conversely, cover the first \(N_\eps\) heights individually and the remaining interval: \(h_{N_\eps}/\eps\asymp N_\eps\). This proves \eqref{eq:height-log}.

The height function is eventually convex, so, after discarding finitely many heights, its gaps are positive, non-increasing and summable. Their tail sums recover the remaining heights. Angelini and Molter's Theorem~4.3 therefore gives \(\dim_\theta V=\theta\beta/[1-(1-\theta)\beta]\) for \(\theta>0\); finite removal and rescaling do not change these dimensions. At \(\theta=0\), use countability.
\end{proof}

For the product core, we use Proposition 2.5 of Falconer, Fraser and Kempton \cite{falconer2020}. Its product bounds determine \(\dim_\theta(X\times V)\) when \(\dimH X=\dimB X\) and \(\dim_\theta V\) exists.

\begin{corollary}[Dimension transitions]\label{cor:power}
For cube or carpet layers at \(h_n=n^{-\alpha}\), put \(D=\max\{d,s+(\alpha+1)^{-1}\}\). Then
\[
\dimB F=D,\qquad \cov_\eps(F)\asymp\eps^{-D},\qquad
|(F)_\eps|\asymp\eps^{d+1-D}.
\]
Thus \(F\) is Minkowski nondegenerate at \(D\), and \eqref{eq:exact-content} applies when \(s+(\alpha+1)^{-1}<d\). For cube layers with either height law of Lemma~\ref{lem:heights},
\begin{equation}\label{eq:intermediate-power}
\dim_\theta F=\max\left\{d,s+\frac\theta{\alpha+\theta}\right\},
\qquad0\le\theta\le1.
\end{equation}
\end{corollary}
\begin{proof}
Combine Lemma~\ref{lem:heights} with Corollary~\ref{cor:prefractal-profile}. For cubes, \(\dimH X=\dimB X=s\), so all intermediate dimensions of \(X\) equal \(s\). The cited product bounds give \(\dim_\theta(X\times V)=s+\dim_\theta V\) for \(\theta>0\); apply Theorem~\ref{thm:summable-approximation}. At \(\theta=0\), use \(\dimH F=d\). All comparisons are as \(\eps\downarrow0\), with constants depending on the construction and height parameters.
\end{proof}

For carpets the exact transfer remains \(\max\{2,\dim_\theta(X\times V)\}\), with separate lower and upper versions. The known carpet spectrum \cite{banaji_kolossvary2024} does not justify substituting \(s_B\) into \eqref{eq:intermediate-power}.

Returning to \(F_\alpha\) from \eqref{eq:cantor-model}, take \(d=1\), \(b=3\), \(M=2\), hence \(s=\log_3 2\) and \(q=2/3\). Its box-dimension gap occurs precisely when
\begin{equation}\label{eq:cantor-threshold}
0<\alpha<\alpha_*:=\frac{s}{1-s}\approx1.70951.
\end{equation}
The product core \(C\times V_\alpha\) has Hausdorff dimension \(s\) and box dimension \(s+(\alpha+1)^{-1}\). It has its own gap for every \(\alpha>0\), whereas \(F_\alpha\) has a gap only when this box dimension exceeds \(1\). The base \(C\times\{0\}\) has no gap. These distinguish the threshold criterion \eqref{eq:product-core-gap} from a comparison of the sets' own dimension gaps.

The link \eqref{eq:density-product-bridge} also computes the moving scaling-core exponents: for every fixed \(\kappa>0\) and \(0<\theta\le1\) such that \(s+\theta/(\alpha+\theta)>1\),
\begin{equation}\label{eq:cantor-density-exponents}
\underline D_{\theta,\kappa}(F_\alpha)
=\overline D_{\theta,\kappa}(F_\alpha)
=s+\frac\theta{\alpha+\theta}.
\end{equation}

For \(\alpha>\alpha_*\), the Minkowski content is \(\mathcal M^1(F_\alpha)=4\). Now replace the heights by the power-logarithmic law of Lemma \ref{lem:heights} at \(\alpha=\alpha_*\). Varying \(\lambda\) leaves all Hausdorff, packing, box and intermediate dimensions equal to \(1\), but changes the normalized tube growth. With \(\beta_*=1-s\), the resulting layered set \(F\) satisfies
\begin{equation}\label{eq:log-tube}
|(F)_\eps|\asymp\eps\bigl[1+(\log(1/\eps))^{\lambda\beta_*}\bigr].
\end{equation}
Thus \(\lowM^1(F)=\infty\) for \(\lambda>0\), \(F\) is Minkowski nondegenerate for \(\lambda=0\), and \(\mathcal M^1(F)=4\) for \(\lambda<0\). Analogous logarithmic corrections occur in nests, chirps and spirals \cite{milicic2018,rolin_vlah_zupanovic2021}.

Finally, geometric heights \(h_n=\rho^n\), \(0<\rho<1\), give \(\cov_\eps(V)\asymp\log(1/\eps)\). An eventual ratio bound \(h_{n+1}\le\rho h_n\) gives the upper estimate \(O(\log(1/\eps))\), with dependence on \(\rho\) and the initial segment. In either horizontal construction, \(s<d\) therefore yields box dimension \(d\) and the content in \eqref{eq:exact-content}. These ratio-lacunary estimates belong to the same sequence theory \cite{elezovic_zupanovic_zubrinic2007,mardesic_resman_zupanovic2012}.
\section*{Acknowledgements}
\label{sec:orgf2b5242}
Part of this work was carried out at the Abdus Salam International Centre for
Theoretical Physics (ICTP). The author acknowledges support through internal project
funding of the Faculty of Informatics, Juraj Dobrila University of Pula.
\section*{Use of AI tools}
\label{sec:orga2ad333}
OpenAI ChatGPT (Codex) assisted with theorem and proof development, analytic checks,
literature review, exposition and figure code. This included the scaling-core,
intermediate-dimension, realization and additive-tube results.

\bibliographystyle{plain}
\bibliography{layered_prefractals}
\end{document}